# On the abundance of the irrational numbers

W. Mückenheim

University of Applied Sciences, Baumgartnerstraße 16, D-86161 Augsburg, Germany

mueckenh@rz.fh-augsburg.de

________________________________

**Abstract.** Based upon the axiom of choice it is proved that the cardinality of the rational numbers is not less than the cardinality of the irrational numbers. This contradicts a main result of transfinite set theory and shows that the axiom of choice is invalid.

## 1. Introduction

In mathematics, as in common life, we find diverging opinions with respect to one or the other question. Some questions are as yet undecided, some questions have been settled or decisions are in progress while others have been proven undecidable. But to mathematics it is unique, that two absolutely contrary opinions do not logically exclude each other but exist simultaneously while there seems to be no chance to pick out a false one and to establish a remaining truth. This case is realised by the philosophy and mathematics of the infinite. Meanwhile it is pending for more than one hundred years. While transfinite set theory is impossible without different degrees of infinity, constructivists and intuitionists deny this notion without running into inconsistencies as is admitted by some of the foremost set theorists:

> ... the attitude of the (neo-)intuitionists that there do not exist altogether non-equivalent infinite sets is consistent, though almost suicidal for mathematics. [1, p. 62]

It would not be astonishing if in different axiomatic systems different results were obtained with respect to peculiarities of those systems. But set theorists on one side and constructivists and intuitionists on the other are certainly believing to address the same entities when speaking of "rational numbers" or of "irrational numbers". In spite of that, the former are convinced that there are infinitely many more irrational numbers than rational numbers while the latter deny that:

> Hence the continua of Weyl, Lebesgue, Lusin, etc. are denumerable ... [2, p. 255]

This situation yields bewildering results:

> Feferman and Levy showed that one cannot prove that there is any non-denumerable set of real numbers which can be well ordered. ... Moreover, they also showed that the statement that the set of all real numbers is the union of a denumerable set of denumerable sets cannot be refuted. [2, p. 62]

> Nevertheless, the great majority of mathematicians refuse to accept the thesis that Cantor's ideas were but a pathological fancy. Though the foundations of set theory are still somewhat shaky, ... [2, p. 345]

Most surprising and by no means to be expected of a pupil of Fraenkel's is that Robinson states:

> Infinite totalities do not exist in any sense of the word (i.e. either really or ideally). More precisely, any mention, or purported mention, of infinite totalities is, literally, *meaningless*. Nevertheless, we should act as if infinite totalities really existed. [3]

Does there exist a correct and an incorrect position? And, if so, who is right, who is wrong?

> Was speziell die beiden für die Mengenlehre besonders wichtigen Axiome "der Auswahl" und "des Unendlichen" betrifft, ... so stellen sich Russell und Whitehead [4] auf den Standpunkt, daß diese Axiome nicht logisch beweisbar sind und (letzten Endes nach subjektivem Ermessen) angenommen oder abgelehnt werden können. Je nachdem man sich für das eine oder das andere entscheidet, gestaltet sich der Bereich der Mathematik enger oder weiter in einem Maß, das vermöge der Anlage des Werks überall sichtbar wird. Man wird so in den Stand gesetzt, den Wert und die Notwendigkeit der Voraussetzungen auf Grund der aus ihnen herleitbaren *Folgerungen* zu beurteilen, ein bei der Grundlegung der Mathematik vielfach durchaus angemessenes und fast unentbehrliches Verfahren.[1] [5, p. 182]

We will follow this advice of Fraenkel, namely to judge about the value and necessity of the basic axioms, in particular of the axiom of choice, by considering its consequences, in order to settle this question. These consequences will turn out to entail what, in an euphemistic way, by set theorists usually is called a "paradoxical result", in order to avoid the term self-contradiction.

> Apart from the well-ordering theorem some statements of quite different character - in particular geometrical statements - have been proved by means of the axiom of choice, which because of their paradoxical character induced some mathematicians to reject the axiom. Presumably the earliest statement (1914) of this kind is Hausdorff's discovery that half of the sphere's surface is congruent to a third of it. ... It may surprise scholars working in the field ... that even after more than half a century of utilising the axiom of choice and well-ordering theorem, a number of first-rate mathematicians (especially French) have not essentially changed their distrustful attitude. [2, p. 83]

---

[1] Concerning the particularly important axioms of set theory, namely the axioms "of choice" and "of infinity", ... Russell and Whitehead [4] take the position that these axioms cannot be proved logically and (after all, by subjective choice) may be accepted or refused. Due to one or the other decision, the field of mathematics becomes narrower or wider, as, by the very concept of this book becomes visible in all of its parts. So one is enabled to judge about value and necessity of the assumptions from their *consequences* - with respect to the foundations of mathematics often an absolutely appropriate and nearly indispensable method.

Transfinite set theory arises from Cantor's observation [6-8] that the set of all irrational numbers has infinitely many more members than the set of all rational numbers. While the latter has the same cardinality $\aleph_0$ as the set $\mathbb{N}$ of all natural numbers n, the cardinality $\aleph$ of the set of all irrational numbers is larger, $\aleph = 2^{\aleph_0}$. It is proven to be uncountable, i.e., any bijection with $\mathbb{N}$ can be excluded.

## 2. A common-sense approach

The rational numbers and the irrational numbers make up the linear continuum. The irrational numbers are commonly understood as Dedekind-cuts, i.e., as the gaps between the rational numbers. The number of gaps, however, cannot be larger than the number of their ends, which, though open intervals are involved, are uniquely defined by the rational numbers. We can deduce this fact "without arbitrariness from the general principles of the laws of logic"[2] or "based upon a common logical principle which is already necessary and indispensable for the basics of mathematical conclusion"[3]. But this fundamental problem is usually not clearly spelled out. It is neither addressed by Fraenkel [5] nor can it be found in most other text books on set theory. One exception is made, at least indirectly, in a later version of Fraenkel's text book edited by Levy:

> Since the set of the irrational numbers has the cardinal $\aleph$, in our case the set of the gaps in [the rational numbers] D has a greater cardinal than D itself; this phenomenon is not surprising since the cuts are essentially subsets and not members of the set. [1, p. 160]

The problem is mentioned though it is immediately relativated because otherwise it would be suicidal for set theory. We note that, as the number of gaps is countable, at least one of those subsets must include uncountably many members, i.e., irrational numbers without one single rational number among them. In the same book we find, only five pages later, the order-theoretic definition of the real numbers C:

> C has a denumerable subset D such that between every two members of C there is a member of D. [1, p. 165]

---
[2] quoting a statement of Zermelo [9, p. 263], originally concerning the question of definiteness.
[3] quoting a statement of Hilbert [10, p. 152], originally aimed at the axiom of choice.

Now we have to assume a rational between any pair of irrational numbers, and no longer between (infinite) subsets only. These are clearly two contradictory statements, one of them being necessarily false. But common sense is not highly esteemed in set theory. By definition, only a bijective mapping can prove equivalence. We will not consider the question here, why a bijective mapping keeps its power in the infinite, while most laws and properties of finite numbers are no longer valid there.

### 3. Constructing a bijection

We will construct a bijection between the set $\mathbb{X}$ of all positive irrational numbers and a subset of the set $\mathbb{Q}$ of all positive rational numbers. (The case of negative rational and irrational numbers can then be treated in a similar way.) The proof is based upon the following assumptions:

1) It is always possible to distinguish between a rational number and an irrational number and between a number with index and a number without index.

2) There exists always a rational number between two real numbers.

3) We accept the axiom of choice, according to which there exists a well-ordering of the complete set $\mathbb{X}$ of all positive irrational numbers (though no prescription is given how to construct it). Further the axiom of choice transforms the infinitely many successive arbitrary steps, required to exhaust $\mathbb{X}$, into simultaneous ones (though no prescription is given, how this is accomplished). These infinitely many steps are needed to biject $\mathbb{X}$ with a set of ordinal numbers $\{\alpha_\nu\}$ of sufficient cardinality $|\{\alpha_\nu\}| \geq |\mathbb{X}|$.

These assumptions are sufficient to prove the

**Theorem.** $|\mathbb{Q}| \geq |\mathbb{X}|$.

**Proof.** Define two sets, one of them, Q, containing only the number zero, Q = {0}, and the other one, X, being empty, X = {}. Let the elements ξ of 𝕏 be well-ordered and add them, maintaining their order, one after the other, to the set X. Every time when one of these ξ is transferred from 𝕏 into X, a rational number q' ∉ Q which always exists between ξ and an also always existing element q ∈ Q, satisfying q < ξ,

$$q < q' < \xi \tag{1}$$

is transferred from ℚ into Q. In the first step, we have obviously q = 0, and q' can be any positive rational number q' < ξ, where ξ is the first element of 𝕏 according to the assumed well ordering. Enumerate q' by the first of the ordinal numbers $\alpha_\nu$ and add it to Q. Enumerate ξ in the same way and add it to X. Then take the next ξ ∈ 𝕏. A not yet enumerated rational number q' ∉ Q exists between ξ and an already enumerated rational number q ∈ Q, satisfying q < ξ, according to (1). Enumerate q' by the next of the ordinal numbers $\alpha_\nu$ and add it to Q. Enumerate ξ in the same way and add it to X.

Continue this process until the set 𝕏 is exhausted. According to the axiom of choice this is possible, unless the sets ℚ and $\{\alpha_\nu\}$ were exhausted prematurely, i.e., if there were no longer a rational number available, serving as a q', or no longer an $\alpha_\nu$ available, serving as an index. But that cannot happen, because there is *always* a rational number q' ∉ Q between any two real numbers q and ξ, because there exist additional rational numbers q' between *any* two members of Q. The set $\{\alpha_\nu\}$ cannot become exhausted prematurely too, because it contains enough members by definition. Therefore, the required rational partner $q_{\alpha_\nu}$ of $\xi_{\alpha_\nu}$ does always exist as well as its index. Though it is clear that an infinite process in fact can never be carried out, the axiom of choice enables us to assume its feasibility and the *existence* of infinitely many pairs $(q_{\alpha_\nu}, \xi_{\alpha_\nu})$.

Finally we have $X_{fin} = \mathbb{X}$ and $Q_{fin} \subseteq \mathbb{Q}$. The bijection results in $|Q_{fin}| = |X_{fin}|$ and, hence, $|\mathbb{Q}| \geq |\mathbb{X}|$.

The definition of the ordinal type of the real numbers, our assumption 2, requires a rational number between any two real numbers and, therefore, it also requires a rational number between any two elements of $Q_{fin}$. As already mentioned, after $\mathbb{X}$ has been exhausted, there must remain enough rational numbers to be placed between any two adjacent elements of $Q_{fin}$, doubling its cardinal number. After the cardinal number of $Q_{fin}$ has been doubled, the very same argument leads to twice the number of elements again and again; this procedure may be applied repeatedly for an infinite number, $\aleph_0$, of times, finally yielding the cardinal number of the set of all positive rational numbers

$$|\mathbb{Q}| \geq 2^{\aleph_0}|Q_{fin}| = 2^{\aleph_0}|\mathbb{X}|. \qquad (2)$$

### 4. Conclusion

This result becomes easily comprehensible if we adopt Dedekind's position, who claims that the irrational numbers have to be created.

> Jedesmal nun, wenn ein Schnitt ($A_1/A_2$) vorliegt, welcher durch keine rationale Zahl hervorgebracht wird, so erschaffen wir eine neue, eine irrationale Zahl ...[4] [11]

Such a procedure could not have yielded more than finitely many irrational numbers during the life time of the universe. Nevertheless, our result is in striking contradiction with transfinite set theory. There remain two alternatives: Either the result is correct and transfinite set theory is false, or at least one of our three assumptions is wrong. The first assumption is usually taken for granted. It is certainly required, if set-theoretic considerations involving sets like $\mathbb{Q}$ and $\mathbb{X}$ in general should supply meaningful results. The second one can easily be verified and is correct without any doubt: Two irrational numbers without a rational number between them are not different and, hence, are not two but only one and the same irrational number. Therefore we have to question our assumption 3, i.e., the axiom of choice. It is responsible for the above proof by facilitating well-ordering and exhaustibility of the set of all positive irrational numbers. It cannot be maintained. In particular its implication of actual existence and simultaneous treatability of infinite sets would allow us, for instance, to take all rational numbers, multiply each one by $e$ or $\pi$ and create as many (and more) irrational numbers as there are rational numbers. This result could again be defeated by the above proof, based on the axiom of choice.

---

[4] Now, every time when there is a cut ($A_1/A_2$) which is not generated by a rational number, then we create a new, an irrational number ...

As Fraenkel asserts, the appearance of paradoxes may be used to positively prove the non-existence of the respective set.

> Vor allem ist die Bildung der wichtigsten Klasse paradoxer Mengen, nämlich der allzu umfassenden Mengen (Antinomien von Burali-Forti, Russell usw.) durch unsere Axiome ausgeschlossen. ... Umgekehrt lassen sich daher die bei gewissen Paradoxien auftretenden Gedankengänge benutzen, um positiv die Nichtexistenz der betreffenden paradoxen Menge nachzuweisen [5]. [5, p. 214]

This positive proof has been given in case of the set of all positive irrational numbers. If the axiom of choice is maintained, we obtain a paradoxical result. If this axiom is dropped, there is no chance of well-ordering the set of all irrational numbers. Hence, we cannot attribute any cardinal number to it; the set of all irrational numbers, and with it the set of all real numbers, is set-theoretically untreatable.

The impossible set of all non-finite numbers, however, may enjoy a revival because it might well turn out to contain only one element, most conveniently denoted by $\infty$.

---

[5] In particular the creation of the most important class of paradox sets, namely the much too large sets (antinomies of Burali-Forti, Russell etc.) is excluded by our axioms. ... Vice versa, we can use the ideas connected with certain paradoxes in order to positively prove the nonexistence of the respective set.